\documentclass[11pt]{article}
\usepackage{amsfonts}
\usepackage{amsfonts}
\usepackage{amsfonts}
\usepackage{mathrsfs}
\usepackage{bm}
\usepackage{amssymb,amsmath} 
 \allowdisplaybreaks

\catcode`!=11
\let\!int\int \def\int{\displaystyle\!int}
\let\!lim\lim \def\lim{\displaystyle\!lim}
\let\!sum\sum \def\sum{\displaystyle\!sum}
\let\!sup\sup \def\sup{\displaystyle\!sup}
\let\!inf\inf \def\inf{\displaystyle\!inf}
\let\!cap\cap \def\cap{\displaystyle\!cap}
\let\!max\max \def\max{\displaystyle\!max}
\let\!min\min \def\min{\displaystyle\!min}
\let\!frac\frac \def\frac{\displaystyle\!frac}
\catcode`!=12

\let\oldsection\section
\renewcommand\section{\setcounter{equation}{0}\oldsection}

\allowdisplaybreaks
\def\pf{\it{Proof.}\rm\quad}

\def\N{\mathbb{N}}\def\Z{\mathbb{Z}}

\newtheorem{thm}{Theorem}[section]
\newtheorem{lem}[thm]{Lemma}
\newtheorem{cor}[thm]{Corollary}

\newtheorem{ex}{Example}[section]
\begin{document}
\title {\bf Integrals of logarithmic functions and alternating multiple zeta values}
\author{
{Ce Xu\thanks{Corresponding author. Email: 15959259051@163.com (C. XU)}}\\[1mm]
\small School of Mathematical Sciences, Xiamen University\\
\small Xiamen
361005, P.R. China}
\date{}
\maketitle \noindent{\bf Abstract }
By using the method of iterated integral representations of series, we establish some explicit relationships between multiple zeta values and Integrals of logarithmic functions. As applications of these relations, we show that multiple zeta values of the form
\[\zeta ( {\bar 1,{{\left\{ 1 \right\}}_{m - 1}},\bar 1,{{\left\{ 1 \right\}}_{k - 1}}} ),\ (k,m\in \N)\]
for $m=1$ or $k=1$, and
\[\zeta ( {\bar 1,{{\left\{ 1 \right\}}_{m - 1}},p,{{\left\{ 1 \right\}}_{k - 1}}}),\ (k,m\in\N)\]
for $p=1$ and $2$, satisfy certain recurrence relations which allow us to write them in terms of zeta values, polylogarithms and $\ln 2$. Moreover, we also prove that the multiple zeta values $\zeta ( {\bar 1,{{\left\{ 1 \right\}}_{m - 1}},3,{{\left\{ 1 \right\}}_{k - 1}}} )$ can be expressed as a rational linear combination of products of zeta values, multiple polylogarithms and $\ln 2$ when $m=k\in \N$. Furthermore, we also obtain reductions for certain multiple polylogarithmic values at $\frac {1}{2}$.
\\[2mm]
\\[2mm]
\section{Introduction}
Let $\N$ be the set of natural numbers. For integers $n,m\in \N$, a generalized harmonic number $\zeta_n(m)$ (also called the partial sums of Riemann zeta function when $m>1$) is defined by (see \cite{Xu2016,X2016})
\[{\zeta _n}\left( m \right): = \sum\limits_{j = 1}^n {\frac{1}{{{j^m}}}}, m\in \N\tag{1.1}\]
which is a natural generalization of the harmonic number \cite{FS1998,F2005}
\[{H_n}: = {\zeta _n}\left( 1 \right)= \sum\limits_{j = 1}^n {\frac{1}{j}} .\tag{1.2}\]
Similarly, let
\[{L_n}\left( m \right) := \sum\limits_{j = 1}^n {\frac{{{{( - 1)}^{j - 1}}}}{{{j^m}}}, m\in \N} \tag{1.3}\]
denote the alternating harmonic numbers (also called the partial sums of alternating Riemann zeta function, see \cite{Xu2016}). Here the Riemann zeta function and alternating zeta function are defined respectively by the convergent series \cite{A2000,B1985,B1989}
 \[\zeta(s):=\sum\limits_{n = 1}^\infty {\frac {1}{n^{s}}},\Re(s)>1, \tag{1.4}\]
and
\[\bar \zeta \left( s \right) = \sum\limits_{n = 1}^\infty  {\frac{{{{\left( { - 1} \right)}^{n - 1}}}}{{{n^s}}}} ,\;{\mathop{\Re}\nolimits} \left( s \right) \ge 1.\tag{1.5}\]
In general, for $m\in \N,\ {\bf S}:=(s_1,s_2,...,s_m)\in (\N)^m$, and a non-negative integer $n$, the multiple harmonic number (MHN for short) is defined by
\[{\zeta _n}\left( {{s_1},{s_2}, \ldots ,{s_m}} \right) := \sum\limits_{1 \le {k_m} <  \cdots  < {k_1} \le n} {\frac{1}{{k_1^{{s_1}} \cdots k_m^{{s_m}}}}}\tag{1.6} \]
and the multiple harmonic star number (MHSN  for short) is given by
\[{\zeta^\star_n}\left( {{s_1},{s_2}, \ldots ,{s_m}} \right) := \sum\limits_{1 \le {k_m} \le  \cdots  \le {k_1} \le n} {\frac{1}{{k_1^{{s_1}} \cdots k_m^{{s_m}}}}}.\tag{1.7}\]
The integers $m$ and $\omega : = \left| {\bf S} \right|: = \sum\nolimits_{i = 1}^m {{s_i}}$ are called the depth and the weight of a multiple harmonic number or multiple harmonic star number. By convention, we put ${\zeta _n}\left( {\bf S} \right)=0$ if $n<m$, and ${\zeta _n}\left(\emptyset \right)={\zeta^\star _n}\left(\emptyset \right)=1$. By ${\left\{ {{s_1},{s_2}, \ldots ,{s_j}} \right\}_m}$ we denote the sequence of depth $mj$ with $m$ repetitions of ${\left( {{s_1},{s_2}, \ldots ,{s_j}} \right)}.$\ (Many papers use the opposite convention, with the $k_i$'s ordered by $n\geq k_m >  \cdots > k_ 1\geq 1$ or $n\geq k_m \geq  \cdots \geq k_ 1\geq 1$, (see \cite{L2012,KP2013,DZ1994,DZ2012})). Moreover, we put a bar on top of $s_j\ (j=1,2,\cdots, m)$ if there is a sign $(-1)^{k_j}$ appearing in the denominator on the right. For example
\[{\zeta _n}\left( {{{\bar s}_1},{s_2}, \ldots ,{{\bar s}_m}} \right) = \sum\limits_{1 \le {k_m} <  \cdots  < {k_1} \le n} {\frac{{{{\left( { - 1} \right)}^{{k_1} + {k_m}}}}}{{k_1^{{s_1}} \cdots k_m^{{s_m}}}}},\tag{1.8}\]
\[{\zeta^\star _n}\left( {{s_1},{{\bar s}_2}, \ldots ,{{\bar s}_m}} \right) = \sum\limits_{1 \le {k_m} \le  \cdots  \le {k_1} \le n} {\frac{{{{\left( { - 1} \right)}^{{k_2} + {k_m}}}}}{{k_1^{{s_1}} \cdots k_m^{{s_m}}}}} .\tag{1.9}\]
The limit cases of MHNs or MHSNs give rise to multiple zeta values (MZVs for short, also called $m$-fold Euler sums [10]) or multiple zeta-star values (MZSVs for short):
\[\zeta \left( {{s_1},{s_2}, \ldots ,{s_m}} \right) := \mathop {\lim }\limits_{n \to \infty } {\zeta _n}\left( {{s_1},{s_2}, \ldots ,{s_m}} \right) = \sum\limits_{{k_1} >  \cdots  > {k_m} \ge 1} {\frac{1}{{k_1^{{s_1}} \cdots k_m^{{s_m}}}}},\tag{1.10}\]
\[\zeta^\star\left( {{s_1},{s_2}, \ldots ,{s_m}} \right) := \mathop {\lim }\limits_{n \to \infty } {\zeta^\star _n}\left( {{s_1},{s_2}, \ldots ,{s_m}} \right) = \sum\limits_{{k_1} \ge  \cdots  \ge {k_m} \ge 1} {\frac{1}{{k_1^{{s_1}} \cdots k_m^{{s_m}}}}}\tag{1.11} \]
defined for $s_2,\ldots, s_m\geq 1$ and $s_1\geq 2$ to ensure convergence of the series. For alternating MHNs or MHSNs, we have the similar limit cases. For example
\[\zeta \left( {{{\bar s}_1},{s_2}, \ldots ,{{\bar s}_m}} \right) := \mathop {\lim }\limits_{n \to \infty } {\zeta _n}\left( {{{\bar s}_1},{s_2}, \ldots ,{{\bar s}_m}} \right) = \sum\limits_{1 \le {k_m} <  \cdots  < {k_1}} {\frac{{{{\left( { - 1} \right)}^{{k_1} + {k_m}}}}}{{k_1^{{s_1}} \cdots k_m^{{s_m}}}}},\tag{1.12} \]
\[{\zeta ^ \star }\left( {{s_1},{{\bar s}_2}, \ldots ,{{\bar s}_m}} \right): = \mathop {\lim }\limits_{n \to \infty } \zeta _n^ \star \left( {{s_1},{{\bar s}_2}, \ldots ,{{\bar s}_m}} \right) =\sum\limits_{1 \le {k_m} \le  \cdots  \le {k_1} \le n} {\frac{{{{\left( { - 1} \right)}^{{k_2} + {k_m}}}}}{{k_1^{{s_1}} \cdots k_m^{{s_m}}}}}.\tag{1.13} \]
A good deal of work on multiple zeta values has focused on the problem of determining when ¡®complicated¡¯ sums can be expressed in terms of ¡®simpler¡¯ sums. A crude but convenient measure of the complexity of the sums (1.10)-(1.13) is the number $m$ of nested summations. This is also equal to the number of arguments in the definitions
(1.10)-(1.13), and is called the depth. Thus, researchers are interested in determining which
sums can be expressed in terms of other sums of lesser depth.
When $m=2$, the multiple zeta values $\zeta ( {{s_1},{s_2}} )$ (or $\zeta^\star( {{s_1},{s_2}} )$) are also called double linear Euler sums \cite{FS1998}. The general double nonlinear Euler sums
involving harmonic number and alternating harmonic number are defined by the series (see \cite{Xu2016,X2016})
$$\sum\limits_{n = 1}^\infty  {\frac{{\prod\limits_{i = 1}^{{m_1}} {\zeta _n^{{q_{_i}}}\left( {{k_i}} \right)\prod\limits_{j = 1}^{{m_2}} {L_n^{{l_j}}\left( {{h_j}} \right)} } }}{{{n^p}}}} ,\ \sum\limits_{n = 1}^\infty  {\frac{{\prod\limits_{i = 1}^{{m_1}} {\zeta _n^{{q_{_i}}}\left( {{k_i}} \right)\prod\limits_{j = 1}^{{m_2}} {L_n^{{l_j}}\left( {{h_j}} \right)} {}} }}{{{n^p}}}}{\left( { - 1} \right)}^{n - 1}.$$
where $p (p>1),m_{1},m_{2},q_{i},k_{i},h_{j},l_{j}$ are positive integers. Many values of double nonlinear Euler sums can be expressed as a rational linear combination of zeta values \cite{BBG1994,BBG1995,FS1998,Xu2016,X2016}.
The study of these Euler sums was started by Euler. After that many different methods, including partial fraction expansions, Eulerian Beta integrals, summation formulas for generalized hypergeometric functions and contour integrals, have been used to evaluate these sums. For details and historical introductions,
please see \cite{BBG1994,BBG1995,BBGP1996,BZB2008,FS1998,F2005,M2014,Xu2016,X2016} and references therein. The evaluation of $\zeta \left( {{s_1},{s_2}} \right)$ (or $\zeta^\star\left( {{s_1},{s_2}} \right)$) in terms of values of Riemann zeta function at positive integers is known when $s_2=1,\ s_1=s_2,\ (s_1,s_2)=(2,4),(4,2)$ or $s_1+s_2$ is odd \cite{BBG1994,BBG1995,FS1998}. When $m=3$, the multiple zeta values $\zeta \left( {s_1,s_2,s_3} \right)$ are investigated \cite{BG1996,FS1998,CM1994}. Markett \cite{CM1994} gave explicit reductions to zeta values for all triple sums $\zeta \left( {s_1,s_2,s_3} \right)$ with $s_1+s_2+s_3\le 6$, and he proved an explicit formula for $\zeta \left( {s,1,1} \right)$ in terms of zeta values. In \cite{BG1996}, Borwein and Girgensohn proved that all $\zeta \left( {s_1,s_2,s_3} \right)$ with $s_1+s_2+s_3$ is even or less than or equal to 10 or $s_1+s_2+s_3=12$ were reducible to zeta values and double sums. The best results to date are due to Jonathan M. Borwein et al \cite{BBBL1997}, D. Zagier \cite{DZ2012} and Kh. Hessami Pilehrood et.al \cite{KP2013}. Zagier proved that the multiple zeta values of the form $\zeta \left( {{{\left\{ 2 \right\}}_a},3,{{\left\{ 2 \right\}}_b}} \right)$ (or $\zeta^\star \left( {{{\left\{ 2 \right\}}_a},3,{{\left\{ 2 \right\}}_b}} \right)$) can be expressed in terms of ordinary zeta values and gave explicit formulas:
 \[\zeta \left( {{{\left\{ 2 \right\}}_a},3,{{\left\{ 2 \right\}}_b}} \right) = 2\sum\limits_{r = 1}^{a + b + 1} {{{\left( { - 1} \right)}^r}\left( {\left( {\begin{array}{*{20}{c}}
   {2r}  \\
   {2b + 2}  \\
\end{array}} \right) - \left( {1 - {2^{ - 2r}}} \right)\left( {\begin{array}{*{20}{c}}
   {2r}  \\
   {2a + 2}  \\
\end{array}} \right)} \right)\zeta \left( {2r + 1} \right)H\left( {a + b + 1 - r} \right)} ,\]
where $H\left( m \right) = \zeta \left( {{{\left\{ 2 \right\}}_m}} \right) = \frac{{{\pi ^{2m}}}}{{\left( {2m + 1} \right)!}}$. In [29], Hessami Pilehrood et al. gave some new binomial identities for multiple harmonic numbers ${\zeta _n}\left( {{{\left\{ 1 \right\}}_a},c,{{\left\{ 1 \right\}}_b}} \right),{\zeta _n}\left( {{{\left\{ 2 \right\}}_a},3,{{\left\{ 2 \right\}}_b}} \right)$ and provided a new proof of Zagier's formula for $\zeta^\star \left( {{{\left\{ 2 \right\}}_a},3,{{\left\{ 2 \right\}}_b}} \right)$. For example:
\[{\zeta _n}\left( {{{\left\{ 2 \right\}}_a},1} \right) = 2\sum\limits_{k = 1}^n {\frac{{\left( {\begin{array}{*{20}{c}}
   n  \\
   k  \\
\end{array}} \right)}}{{{k^{2a + 1}}\left( {\begin{array}{*{20}{c}}
   {n + k}  \\
   k  \\
\end{array}} \right)}}},\]
letting $n$ tend to infinity in above equation, then have
\[\zeta \left( {{{\left\{ 2 \right\}}_a},1} \right) = 2\zeta \left( {2a + 1} \right),\ a \in \N.\]
The purpose of the present paper is to study multiple zeta values of the form
\[\zeta \left( {\bar 1,{{\left\{ 1 \right\}}_{m - 1}},\bar 1,{{\left\{ 1 \right\}}_{k - 1}}} \right),\zeta \left( {\bar 1,{{\left\{ 1 \right\}}_{m - 1}},p,{{\left\{ 1 \right\}}_{k - 1}}} \right)\]
for $p=1,2,3$ and $m,k\in \N$.
In this paper, we prove that the multiple zeta values \[\zeta \left( {\bar 1,{{\left\{ 1 \right\}}_{m - 1}},\bar 1} \right),\zeta \left( {\bar 1,\bar 1,{{\left\{ 1 \right\}}_{k - 1}}} \right),\zeta \left( {\bar 1,{{\left\{ 1 \right\}}_k}} \right),\zeta \left( {\bar 1,{{\left\{ 1 \right\}}_{m - 1}},2,{{\left\{ 1 \right\}}_{k - 1}}} \right)\]
 can be expressed as a rational linear combination of products of zeta values, polylogarithms and ln2. For instance,
\[\zeta \left( {\bar 1,{{\left\{ 1 \right\}}_m}} \right) = {\left( { - 1} \right)^{m + 1}}\frac{{{{\ln }^{m + 1}}2}}{{\left( {m + 1} \right)!}},\ \zeta \left( {\bar 1,\bar 1,{{\left\{ 1 \right\}}_m}} \right) =  - {\rm Li}{_{m + 2}}\left( {\frac{1}{2}} \right).\]
The polylogarithm function is defined for $\left| x \right| \le 1$ by
\[{\rm Li}{_p}\left( x \right) = \sum\limits_{n = 1}^\infty  {\frac{{{x^n}}}{{{n^p}}}}, \Re(p)>1,\tag{1.14}\]
with ${\rm Li}_1(x)=-\ln(1-x),\ x\in [-1,1)$. The generalized multiple polylogarithm function and multiple polylogarithm star function are defined by
\[{\rm{L}}{{\rm{i}}_{{s_1},{s_2}, \cdots ,{s_m}}}\left( x \right): = \sum\limits_{1 \le {k_m} <  \cdots  < {k_1}} {\frac{{{x^{{k_1}}}}}{{k_1^{{s_1}}k_2^{{s_2}} \cdots k_m^{{s_m}}}}} ,\;x \in \left[ { - 1,1} \right),\tag{1.15}\]
\[{\rm{Li}}_{{s_1},{s_2}, \cdots ,{s_m}}^ \star \left( x \right): = \sum\limits_{1 \le {k_m} \le  \cdots  \le {k_1}} {\frac{{{x^{{k_1}}}}}{{k_1^{{s_1}}k_2^{{s_2}} \cdots k_m^{{s_m}}}}} ,\;x \in \left[ { - 1,1} \right).\tag{1.16}\]
For convenience, when ${\bf S}:=(s_1,s_2,\ldots,s_m)\in (\Z)^m$ in above definitions of (1.15) and (1.16), we use the following notations
\[\zeta \left( {{s_1},{s_2}, \cdots ,{s_m};x} \right): = {\rm{L}}{{\rm{i}}_{{s_1},{s_2}, \cdots ,{s_m}}}\left( x \right),\;{\zeta ^ \star }\left( {{s_1},{s_2}, \cdots ,{s_m};x} \right): = {\rm{Li}}_{{s_1},{s_2}, \cdots ,{s_m}}^ \star \left( x \right).\tag{1.17}\]
Of course, if $s_1>1$, then we can allow $x=1$. To avoid confusion with the notion of analytic continuation, we shall henceforth adopt the notation of \cite{BBBL1997}, in which each $s_j$ in (1.17) is replaced by ${\overline {-s_j}}$ when $s_j<0$. Thus, for example, $\zeta \left( {\bar 1} \right) =  - \ln 2$.
Furthermore, we show that the multiple zeta values $\zeta \left( {{\overline {k+1}},{{\left\{ 1 \right\}}_{m - 1}}} \right)$ and $\zeta \left( {\bar 1,{{\left\{ 1 \right\}}_{m - 1}},\bar 1,{{\left\{ 1 \right\}}_{k - 1}}} \right)$
can be expressed as a rational linear combination of zeta values, polylogarithms, $\ln2$ and multiple polylogarithmic values at $1/2$. Now we state our main theorems.
\section{Main Theorems}
The main result of this paper can be stated as follows.
\begin{thm} For integers $k\geq 0$ and $m\geq 1$, then the following identity holds:
\[\zeta \left( {\overline {k + 2},{{\left\{ 1 \right\}}_{m - 1}}} \right) = \frac{{{{\left( { - 1} \right)}^{m + k}}}}{{m!k!}}I\left( {k,m} \right),\tag{2.1}\]
where $I(k,m)$ is defined by the integral
\[I\left( {k,m} \right): = \int\limits_0^1 {\frac{{{{\left( {\ln x} \right)}^k}{{\ln }^m}\left( {1 + x} \right)}}{x}dx} .\]
Moreover, we have the following result of the integral $I(k,m)$
\begin{align*}
I\left( {k,m} \right) =& \frac{1}{{m + k + 1}}{\left( {\ln 2} \right)^{m + k + 1}} + \left( {m + k} \right)!\zeta \left( {m + k + 1} \right) \\&- \sum\limits_{l = 0}^{m + k} {l!\left( {\begin{array}{*{20}{c}}
   {m + k}  \\
   l  \\
\end{array}} \right){{\left( {\ln 2} \right)}^{m + k - l}}{\rm{L}}{{\rm{i}}_{l + 1}}\left( {\frac{1}{2}} \right)} \\
& + {\left( { - 1} \right)^k}\sum\limits_{j = 1}^k {j!\left( {k + m - j} \right)!\left( {\begin{array}{*{20}{c}}
   k  \\
   j  \\
\end{array}} \right)} \left\{ \begin{array}{l}
 \zeta \left( {k + m + 2 - j,{{\left\{ 1 \right\}}_{j - 1}}} \right) \\
 \;\;\;\;\; + \zeta \left( {k + m + 1 - j,{{\left\{ 1 \right\}}_j}} \right) \\
 \end{array} \right\}\\
 & - {\left( { - 1} \right)^k}\sum\limits_{j = 1}^k {\sum\limits_{l = 0}^{k + m - j} {j!l!\left( {\begin{array}{*{20}{c}}
   k  \\
   j  \\
\end{array}} \right)\left( {\begin{array}{*{20}{c}}
   {k + m - j}  \\
   l  \\
\end{array}} \right)} } {\left( {\ln 2} \right)^{m + k - j - l}}\left\{ \begin{array}{l}
 \zeta \left( {l + 2,{{\left\{ 1 \right\}}_{j - 1}};\frac{1}{2}} \right) \\
 \;\;\; + \zeta \left( {l + 1,{{\left\{ 1 \right\}}_j};\frac{1}{2}} \right) \\
 \end{array} \right\}.\tag{2.2}
\end{align*}
\end{thm}
\begin{thm} For positive integers $m$ and $k$, then the following identity holds:
\[\zeta \left( {\bar 1,{{\left\{ 1 \right\}}_{m - 1}},\bar 1,{{\left\{ 1 \right\}}_{k - 1}}} \right) = \frac{{{{\left( { - 1} \right)}^{k - 1}}}}{{k!\left( {m - 1} \right)!}}J\left( {k,m - 1} \right) - \sum\limits_{i = 1}^{m - 1} {\frac{{{{\left( {\ln 2} \right)}^i}}}{{i!}}\zeta \left( {\bar 1,{{\left\{ 1 \right\}}_{m - 1 - i}},\bar 1,{{\left\{ 1 \right\}}_{k - 1}}} \right)},\tag{2.3} \]
where $J\left( {k,m} \right)$ is defined by
\[J\left( {k,m} \right): = \int\limits_0^1 {\frac{{{{\ln }^k}\left( {1 - t} \right){{\ln }^m}\left( {1 + t} \right)}}{{1 + t}}dt}.\]
We have the following explicit formulas
\begin{align*}
J\left( {k,m} \right) =& \frac{1}{{m + 1}}{\left( {\ln 2} \right)^{m + k + 1}} + \sum\limits_{i = 1}^k {\sum\limits_{j = 0}^m {{{\left( { - 1} \right)}^{i + j}}i!j!{{\left( {\ln 2} \right)}^{m + k - i - j}}\left( {\begin{array}{*{20}{c}}
   k  \\
   i  \\
\end{array}} \right)\left( {\begin{array}{*{20}{c}}
   m  \\
   j  \\
\end{array}} \right)\zeta \left( {j + 2,{{\left\{ 1 \right\}}_{i - 1}}} \right)} } \\
&- \sum\limits_{i = 1}^k {\sum\limits_{j = 0}^m {\sum\limits_{l = 0}^j {{{\left( { - 1} \right)}^{i + j}}i!l!{{\left( {\ln 2} \right)}^{m + k - i - l}}\left( {\begin{array}{*{20}{c}}
   k  \\
   i  \\
\end{array}} \right)\left( {\begin{array}{*{20}{c}}
   m  \\
   j  \\
\end{array}} \right)\left( {\begin{array}{*{20}{c}}
   j  \\
   l  \\
\end{array}} \right)\zeta \left( {l + 2,{{\left\{ 1 \right\}}_{i - 1}};\frac{1}{2}} \right)} } },\tag{2.4}\\
J\left( {k,m} \right) =& \sum\limits_{i = 0}^k {\sum\limits_{j = 0}^m {\sum\limits_{l = 0}^j {{{\left( { - 1} \right)}^{i + j}}j!l!{{\left( {\ln 2} \right)}^{m + k - j - l}}\left( {\begin{array}{*{20}{c}}
   k  \\
   i  \\
\end{array}} \right)\left( {\begin{array}{*{20}{c}}
   i  \\
   l  \\
\end{array}} \right)\left( {\begin{array}{*{20}{c}}
   m  \\
   j  \\
\end{array}} \right)\zeta \left( {l + 1,{{\left\{ 1 \right\}}_j};\frac{1}{2}} \right)} } } .\tag{2.5}
\end{align*}
\end{thm}
\begin{thm} For integers $m,k\in \N$, then the following identity holds:
\begin{align*}
\zeta \left( {\bar 1,{{\left\{ 1 \right\}}_{m - 1}},2,{{\left\{ 1 \right\}}_{k - 1}}} \right) =& {\left( { - 1} \right)^m}\left( {\begin{array}{*{20}{c}}
   {m + k}  \\
   k  \\
\end{array}} \right)\zeta \left( {\bar 2,{{\left\{ 1 \right\}}_{m + k - 1}}} \right) - \frac{{{{\left( {\ln 2} \right)}^m}}}{{m!}}\zeta \left( {\bar 2,{{\left\{ 1 \right\}}_{k - 1}}} \right)\\
& - \sum\limits_{i = 1}^{m - 1} {\frac{{{{\left( {\ln 2} \right)}^i}}}{{i!}}\zeta \left( {\bar 1,{{\left\{ 1 \right\}}_{m - i - 1}},2,{{\left\{ 1 \right\}}_{k - 1}}} \right)} .\tag{2.6}
\end{align*}
\end{thm}
\begin{thm} For positive integers $m$ and $k$, we have
\begin{align*}
&{\left( { - 1} \right)^m}\zeta \left( {\bar 1,{{\left\{ 1 \right\}}_{m - 1}},3,{{\left\{ 1 \right\}}_{k - 1}}} \right) + {\left( { - 1} \right)^k}\zeta \left( {\bar 1,{{\left\{ 1 \right\}}_{k - 1}},3,{{\left\{ 1 \right\}}_{m - 1}}} \right)\\
& = \frac{{{{\left( { - 1} \right)}^{m + 1}}}}{{m!}}{\left( {\ln 2} \right)^m}\zeta \left( {\bar 3,{{\left\{ 1 \right\}}_{k - 1}}} \right) + \frac{{{{\left( { - 1} \right)}^{k + 1}}}}{{k!}}{\left( {\ln 2} \right)^k}\zeta \left( {\bar 3,{{\left\{ 1 \right\}}_{m - 1}}} \right) + \zeta \left( {\bar 2,{{\left\{ 1 \right\}}_{m - 1}}} \right)\zeta \left( {\bar 2,{{\left\{ 1 \right\}}_{k - 1}}} \right)\\
& + {\left( { - 1} \right)^{m+1}}\sum\limits_{i = 1}^{m - 1} {\frac{{{{\left( {\ln 2} \right)}^i}}}{{i!}}\zeta \left( {\bar 1,{{\left\{ 1 \right\}}_{m - i - 1}},3,{{\left\{ 1 \right\}}_{k - 1}}} \right)}  + {\left( { - 1} \right)^{k+1}}\sum\limits_{i = 1}^{k - 1} {\frac{{{{\left( {\ln 2} \right)}^i}}}{{i!}}\zeta \left( {\bar 1,{{\left\{ 1 \right\}}_{k - i - 1}},3,{{\left\{ 1 \right\}}_{m - 1}}} \right)}.\tag{2.7}
\end{align*}
\end{thm}
We prove Theorem 2.1 in section 3, Theorem 2.2 in section 4, Theorem 2.3 and Theorem 2.4 in section 5.  We will prove the Theorem 2.1-2.4 by the method of iterated integral representations of series.
\section{Proof of Theorem 2.1}
We now prove our Theorems. First, we need to prove the following lemmas. It will be useful in the development of the main theorems.
\begin{lem}
For integer $k>0$ and $x\in [-1,1)$, we have that
\[{\ln ^k}\left( {1 - x} \right) = {\left( { - 1} \right)^k}k!\sum\limits_{n = 1}^\infty  {\frac{{{x^n}}}{n}{\zeta _{n - 1}}\left( {{{\left\{ 1 \right\}}_{k - 1}}} \right)},\tag{3.1} \]
\[s\left( {n,k} \right) = \left( {n - 1} \right)!{\zeta _{n - 1}}\left( {{{\left\{ 1 \right\}}_{k - 1}}} \right).\tag{3.2}\]
where ${s\left( {n,k} \right)}$ is called (unsigned) Stirling number of the first kind (see \cite{L1974}).
\begin{align*}
& s\left( {n,1} \right) = \left( {n - 1} \right)!,\\
&s\left( {n,2} \right) = \left( {n - 1} \right)!{H_{n - 1}},\\
&s\left( {n,3} \right) = \frac{{\left( {n - 1} \right)!}}{2}\left[ {H_{n - 1}^2 - {\zeta _{n - 1}}\left( 2 \right)} \right],\\
&s\left( {n,4} \right) = \frac{{\left( {n - 1} \right)!}}{6}\left[ {H_{n - 1}^3 - 3{H_{n - 1}}{\zeta _{n - 1}}\left( 2 \right) + 2{\zeta _{n - 1}}\left( 3 \right)} \right], \\
&s\left( {n,5} \right) = \frac{{\left( {n - 1} \right)!}}{{24}}\left[ {H_{n - 1}^4 - 6{\zeta _{n - 1}}\left( 4 \right) - 6H_{n - 1}^2{\zeta _{n - 1}}\left( 2 \right) + 3\zeta _{n - 1}^2\left( 2 \right) + 8H_{n - 1}^{}{\zeta _{n - 1}}\left( 3 \right)} \right].
\end{align*}
The Stirling numbers ${s\left( {n,k} \right)}$ of the first kind satisfy a recurrence relation in the form
\[s\left( {n,k} \right) = s\left( {n - 1,k - 1} \right) + \left( {n - 1} \right)s\left( {n - 1,k} \right),\;\;n,k \in \N,\]
with $s\left( {n,k} \right) = 0,n < k,s\left( {n,0} \right) = s\left( {0,k} \right) = 0,s\left( {0,0} \right) = 1$.
\end{lem}
\pf To prove the first identity we proceed by induction on $k$. Obviously, it is valid for $k=1$. For $k>1$ we use the equality
\[{\ln ^{k{\rm{ + }}1}}\left( {1 - x} \right){\rm{ = }} - \left( {k + 1} \right)\int\limits_0^x {\frac{{{{\ln }^k}\left( {1 - t} \right)}}{{1 - t}}dt} \]
and apply the induction hypothesis, by using Cauchy product of power series, we arrive at
\begin{align*}
{\ln ^{k{\rm{ + }}1}}\left( {1 - x} \right){\rm{ = }}& - \left( {k + 1} \right)\int\limits_0^x {\frac{{{{\ln }^k}\left( {1 - t} \right)}}{{1 - t}}dt} \\
& = {\left( { - 1} \right)^{k + 1}}\left( {k + 1} \right)!\sum\limits_{n = 1}^\infty  {\frac{1}{{n + 1}}\sum\limits_{i = 1}^n {\frac{{{\zeta _{i - 1}}\left( {{{\left\{ 1 \right\}}_{k - 1}}} \right)}}{i}} } {x^{n + 1}}\\
& = {\left( { - 1} \right)^{k + 1}}\left( {k + 1} \right)!\sum\limits_{n = 1}^\infty  {\frac{{{\zeta _n}\left( {{{\left\{ 1 \right\}}_k}} \right)}}{{n + 1}}} {x^{n + 1}}.
\end{align*}
Nothing that ${\zeta _n}\left( {{{\left\{ 1 \right\}}_k}} \right) = 0$ when $n<k$. By simple calculation,
we can deduce (3.1). To prove the second identity of our lemma, we use the following equation (\cite{L1974})
\[{\ln ^k}\left( {1 - x} \right) = {\left( { - 1} \right)^k}k!\sum\limits_{n = k}^\infty  {\frac{{s\left( {n,k} \right)}}{{n!}}{x^n}} ,\: - 1 \le x < 1.\tag{3.3}\]
Thus, by comparing the coefficients of $x^n$ in (3.1) and (3.3), we obtain formula (3.2). The proof of lemma 3.1 is thus completed.\hfill$\square$
\begin{lem} For integers $m\geq 0$ and $n\geq1$, then the following integral identity holds:
\[\int\limits_0^x {{t^{n - 1}}{{\left( {\ln t} \right)}^m}} dt = \sum\limits_{l = 0}^m {l!\left( {\begin{array}{*{20}{c}}
   m  \\
   l  \\
\end{array}} \right)\frac{{{{\left( { - 1} \right)}^l}}}{{{n^{l + 1}}}}{{\left( {\ln x} \right)}^{m - l}}{x^n}},\ x\in (0,1).\tag{3.4} \]
\end{lem}
\pf The lemma is almost obvious. By using integration by parts, we may easily deduce the result.\hfill$\square$\\
Putting $x=1$ and $1/2$ in (3.4), we obtain
\begin{align*}
&\int\limits_0^1 {{t^{n - 1}}{{\left( {\ln t} \right)}^m}} dt = m!\frac{{{{\left( { - 1} \right)}^m}}}{{{n^{m + 1}}}},\tag{3.5}\\
&\int\limits_0^{1/2} {{t^{n - 1}}{{\left( {\ln t} \right)}^m}} dt = {\left( { - 1} \right)^m}\sum\limits_{l = 0}^m {l!\left( {\begin{array}{*{20}{c}}
   m  \\
   l  \\
\end{array}} \right){{\left( {\ln 2} \right)}^{m - l}}\frac{1}{{{2^n}{n^{l + 1}}}}} ,\tag{3.6}\\
&\int\limits_{1/2}^1 {{t^{n - 1}}{{\left( {\ln t} \right)}^m}} dt = m!\frac{{{{\left( { - 1} \right)}^m}}}{{{n^{m + 1}}}} + {\left( { - 1} \right)^{m - 1}}\sum\limits_{l = 0}^m {l!\left( {\begin{array}{*{20}{c}}
   m  \\
   l  \\
\end{array}} \right){{\left( {\ln 2} \right)}^{m - l}}\frac{1}{{{2^n}{n^{l + 1}}}}} .\tag{3.7}
\end{align*}
{\bf Proof of Theorem 2.1}. For integers $m\geq 1$ and $k\geq 0$, by using formula (3.1) and (3.5), we have
\begin{align*}
I\left( {k,m} \right) =& \int\limits_0^1 {\frac{{{{\left( {\ln x} \right)}^k}{{\ln }^m}\left( {1 + x} \right)}}{x}dx} \\
 = &{\left( { - 1} \right)^m}m!\sum\limits_{n = 1}^\infty  {\frac{{{\zeta _{n - 1}}\left( {{{\left\{ 1 \right\}}_{m - 1}}} \right)}}{n}{{\left( { - 1} \right)}^n}\int\limits_0^1 {{x^{n - 1}}{{\ln }^k}x} } dx\\
 =& {\left( { - 1} \right)^{m + k}}k!m!\sum\limits_{n = 1}^\infty  {\frac{{{\zeta _{n - 1}}\left( {{{\left\{ 1 \right\}}_{m - 1}}} \right)}}{{{n^{k + 2}}}}{{\left( { - 1} \right)}^n}} \\
  =& {\left( { - 1} \right)^{m + k}}k!m!\zeta \left( {\overline {k + 2},{{\left\{ 1 \right\}}_{m - 1}}} \right).\tag{3.8}
\end{align*}
Hence, by a direct calculation, we obtain (2.1).
To prove (2.2), applying the change of variable $x = {u^{ - 1}} - 1$ to the above integral, we get the identity
\begin{align*}
I\left( {k,m} \right) &= \int\limits_{1/2}^1 {\frac{{{{\ln }^k}\left( {\frac{{1 - u}}{u}} \right){{\ln }^m}\left( {\frac{1}{u}} \right)}}{{u\left( {1 - u} \right)}}du} \\
& = {\left( { - 1} \right)^m}\int\limits_{1/2}^1 {\frac{{{{\left( {\ln \left( {1 - u} \right) - \ln u} \right)}^k}{{\ln }^m}\left( u \right)}}{{u\left( {1 - u} \right)}}du} \\
& = {\left( { - 1} \right)^m}\sum\limits_{j = 0}^k {{{\left( { - 1} \right)}^{k - j}}\left( {\begin{array}{*{20}{c}}
   k  \\
   j  \\
\end{array}} \right)\int\limits_{1/2}^1 {\frac{{{{\ln }^j}\left( {1 - u} \right){{\ln }^{m + k - j}}\left( u \right)}}{{u\left( {1 - u} \right)}}du} } \\
& = {\left( { - 1} \right)^m}\sum\limits_{j = 0}^k {{{\left( { - 1} \right)}^{k - j}}\left( {\begin{array}{*{20}{c}}
   k  \\
   j  \\
\end{array}} \right)\int\limits_{1/2}^1 {\left\{ {\frac{{{{\ln }^j}\left( {1 - u} \right){{\ln }^{m + k - j}}\left( u \right)}}{u} + \frac{{{{\ln }^j}\left( {1 - u} \right){{\ln }^{m + k - j}}\left( u \right)}}{{1 - u}}} \right\}du} } .\tag{3.9}
\end{align*}
Substituting (3.1) and (3.7) into (3.9) yields the desired result. We complete the proof of Theorem 2.1. \hfill$\square$\\
Putting $x=-1$ in (3.1), we get
\[\zeta \left( {\bar 1,{{\left\{ 1 \right\}}_{k-1}}} \right) = {\left( { - 1} \right)^{k}}\frac{{{{\ln }^{k}}2}}{{{k} !}}.\]
Taking $k=0$ and $1$ in (2.1) and (2.2), we can give the following corollaries.
\begin{cor} For integer $m\geq 1$, we have
\begin{align*}
\zeta \left( {\bar 2,{{\left\{ 1 \right\}}_{m-1}}} \right) =& \frac{{{{\left( { - 1} \right)}^{m }}}}{{\left( {m + 1} \right)!}}{\ln ^{m + 1}}2 + {\left( { - 1} \right)^{m }}\left( {\zeta \left( {m + 1} \right) - {\rm{L}}{{\rm{i}}_{m + 1}}\left( {\frac{1}{2}} \right)} \right) \\&- {\left( { - 1} \right)^{m }}\sum\limits_{j = 1}^{m } {\frac{{{{\left( {\ln 2} \right)}^{m + 1 - j}}}}{{\left( {m + 1- j} \right)!}}} {\rm{L}}{{\rm{i}}_j}\left( {\frac{1}{2}} \right).\tag{3.10}
\end{align*}
\end{cor}
\begin{cor} For integer $m\geq1$, we obtain
\begin{align*}
\zeta \left( {\bar 3,{{\left\{ 1 \right\}}_{m - 1}}} \right)
& = {\left( { - 1} \right)^{m + 1}}\left\{ {m\zeta \left( {m + 2} \right) - \frac{{{{\left( {m + 1} \right)}^2}}}{{\left( {m + 2} \right)!}}{{\left( {\ln 2} \right)}^{m + 2}} - \zeta \left( {m + 1,1} \right)} \right\}\\
& \quad+ \frac{{{{\left( { - 1} \right)}^{m + 1}}}}{{m!}}\sum\limits_{l = 0}^m {l!\left( {\begin{array}{*{20}{c}}
   m  \\
   l  \\
\end{array}} \right)} {\left( {\ln 2} \right)^{m - l}}\zeta \left( {l + 1,1;\frac{1}{2}} \right)\\
& \quad- \frac{{{{\left( { - 1} \right)}^{m + 1}}}}{{\left( {m - 1} \right)!}}\sum\limits_{l = 0}^m {l!\left( {\begin{array}{*{20}{c}}
   m  \\
   l  \\
\end{array}} \right)} {\left( {\ln 2} \right)^{m - l}}{\rm{L}}{{\rm{i}}_{l + 2}}\left( {\frac{1}{2}} \right).\tag{3.11}
\end{align*}
\end{cor}
\section{Proof of Theorem 2.2}
{\bf Proof of Theorem 2.2}. To prove the first identity (2.3), we consider the following multiple integral
\[{M_m}\left( k \right): = \int\limits_0^1 {\frac{1}{{1 + {t_1}}}d{t_1}}  \cdots \int\limits_0^{{t_{m - 1}}} {\frac{1}{{1 + {t_m}}}d{t_1}} \int\limits_0^{{t_m}} {\frac{{{{\ln }^k}\left( {1 - {t_{m + 1}}} \right)}}{{1 + {t_{m + 1}}}}} d{t_{m + 1}}.\tag{4.1}\]
By using power series expansion and formula (3.1), we deduce that
\begin{align*}
{M_m}\left( k \right) =& {\left( { - 1} \right)^{k + m + 1}}k!\sum\limits_{{n_1},{n_2}, \cdots ,{n_{m + 1}} = 1}^\infty  {{{\left( { - 1} \right)}^{{n_1} +  \cdots  + {n_{m + 1}}}}\frac{{{\zeta _{{n_1} - 1}}\left( {\bar 1,{{\left\{ 1 \right\}}_{k - 1}}} \right)}}{{{n_1}}}} \\
& \times \int\limits_0^1 {t_1^{{n_{m + 1}} - 1}d{t_1}}  \cdots \int\limits_0^{{t_{m - 1}}} {t_m^{{n_1} + {n_2} - 1}d{t_m}} \\
 =& {\left( { - 1} \right)^{k + m + 1}}k!\sum\limits_{{n_1},{n_2}, \cdots ,{n_{m + 1}} = 1}^\infty  {{{\left( { - 1} \right)}^{{n_1} +  \cdots  + {n_{m + 1}}}}\frac{{{\zeta _{{n_1} - 1}}\left( {\bar 1,{{\left\{ 1 \right\}}_{k - 1}}} \right)}}{{{n_1}\left( {{n_1} + {n_2}} \right) \cdots \left( {{n_1} +  \cdots  + {n_{m + 1}}} \right)}}} \\
 =& {\left( { - 1} \right)^{k + m + 1}}k!\sum\limits_{{n_1} >  \cdots  > {n_{m + 1}} \ge 1}^\infty  {\frac{{{\zeta _{{n_{m + 1}} - 1}}\left( {\bar 1,{{\left\{ 1 \right\}}_{k - 1}}} \right)}}{{{n_1}{n_2} \cdots {n_{m + 1}}}}{{\left( { - 1} \right)}^{{n_1}}}} \\
 =& {\left( { - 1} \right)^{k + m + 1}}k!\zeta \left( {\bar 1,{{\left\{ 1 \right\}}_m},\bar 1,{{\left\{ 1 \right\}}_{k - 1}}} \right).\tag{4.2}
\end{align*}
Hence, ${M_{m - 1}}\left( k \right) = {\left( { - 1} \right)^{k + m}}k!\zeta \left( {\bar 1,{{\left\{ 1 \right\}}_{m - 1}},\bar 1,{{\left\{ 1 \right\}}_{k - 1}}} \right)$.
On the other hand, by using integration by parts, we have
\begin{align*}
{M_{m - 1}}\left( k \right) =& \ln 2{M_{m - 2}}\left( k \right) - \int\limits_{0 < {t_{m - 1}} <  \cdots  < {t_1} < 1} {\frac{{\ln \left( {1 + {t_1}} \right){{\ln }^k}\left( {1 - {t_{m - 1}}} \right)}}{{\left( {1 + {t_1}} \right)\left( {1 + {t_2}} \right) \cdots \left( {1 + {t_{m - 1}}} \right)}}d{t_1} \cdots d{t_{m - 1}}} \\
 = &\ln 2{M_{m - 2}}\left( k \right) - \frac{1}{2}{\ln ^2}2{M_{m - 3}}\left( k \right)\\
 & - \frac{1}{2}\int\limits_{0 < {t_{m - 2}} <  \cdots  < {t_1} < 1} {\frac{{{{\ln }^2}\left( {1 + {t_1}} \right){{\ln }^k}\left( {1 - {t_{m - 1}}} \right)}}{{\left( {1 + {t_1}} \right)\left( {1 + {t_2}} \right) \cdots \left( {1 + {t_{m - 2}}} \right)}}d{t_1} \cdots d{t_{m - 2}}} \\
 =&\cdots\\
 =&\sum\limits_{i = 1}^{m - 1} {{{\left( { - 1} \right)}^{i - 1}}\frac{{{{\left( {\ln 2} \right)}^i}}}{{i!}}{M_{m - i - 1}}\left( k \right)}  + \frac{{{{\left( { - 1} \right)}^{m - 1}}}}{{\left( {m - 1} \right)!}}\int\limits_0^1 {\frac{{{{\ln }^k}\left( {1 - {t_1}} \right){{\ln }^{m - 1}}\left( {1 + {t_1}} \right)}}{{1 + {t_1}}}d{t_1}} .\tag{4.3}
\end{align*}
Combining (4.2) with (4.3) yields the desired result. To prove the second identity (2.4), applying the change of variable $t = 2u - 1$ to the integral $J\left( {k,m} \right)$, then we easily obtain
\begin{align*}
J\left( {k,m} \right) =& \int\limits_{1/2}^1 {\frac{{{{\ln }^k}\left( {2 - 2u} \right){{\ln }^m}\left( {2u} \right)}}{u}du} \\
 =&\sum\limits_{i = 1}^k {\sum\limits_{j = 0}^m {{{\left( {\ln 2} \right)}^{k + m - i - j}}\left( {\begin{array}{*{20}{c}}
   k  \\
   i  \\
\end{array}} \right)\left( {\begin{array}{*{20}{c}}
   m  \\
   j  \\
\end{array}} \right)} } \int\limits_{1/2}^1 {\frac{{{{\ln }^i}\left( {1 - u} \right){{\ln }^j}\left( u \right)}}{u}du}  \\&+\sum\limits_{j = 0}^m {{{\left( {\ln 2} \right)}^{k + m - j}}\left( {\begin{array}{*{20}{c}}
   m  \\
   j  \\
\end{array}} \right)\int\limits_{1/2}^1 {\frac{{{{\ln }^j}\left( u \right)}}{u}du} } .\tag{4.4}
\end{align*}
Substituting (3.1) and (3.7) into (4.4), we obtain the result. Similarly, to prove the third identity (2.5), applying the change of variable $t = 1-2x$ to the integral $J\left( {k,m} \right)$, we have that
\begin{align*}
J\left( {k,m} \right) &= \int\limits_0^{1/2} {\frac{{{{\ln }^k}\left( {2x} \right){{\ln }^m}\left( {2 - 2x} \right)}}{{1 - x}}dx} \\
& = \sum\limits_{i = 0}^k {\sum\limits_{j = 0}^m {{{\left( {\ln 2} \right)}^{m + k - i - j}}\left( {\begin{array}{*{20}{c}}
   k  \\
   i  \\
\end{array}} \right)} } \left( {\begin{array}{*{20}{c}}
   m  \\
   j  \\
\end{array}} \right)\int\limits_0^{1/2} {\frac{{{{\ln }^i}\left( x \right){{\ln }^j}\left( {1 - x} \right)}}{{1 - x}}dx} .\tag{4.5}
\end{align*}
Combining (3.1), (3.7) and (4.5), we deduce (2.5). The proof of Theorem 2.2 is thus completed.\hfill$\square$\\
Letting $m=0$ in (2.4), we get the recurrence relation
\begin{align*}
&\sum\limits_{i = 1}^k {{{\left( { - 1} \right)}^i}i!{{\left( {\ln 2} \right)}^{k - i}}\left( {\begin{array}{*{20}{c}}
   k  \\
   i  \\
\end{array}} \right)\left\{ {\zeta \left( {2,{{\left\{ 1 \right\}}_{i - 1}}} \right) - \zeta \left( {2,{{\left\{ 1 \right\}}_{i - 1}};\frac{1}{2}} \right)} \right\}} \\
& = {\left( { - 1} \right)^k}k!{\rm{L}}{{\rm{i}}_{k + 1}}\left( {\frac{1}{2}} \right) - {\left( {\ln 2} \right)^{k + 1}}.\tag{4.6}
\end{align*}
On the other hand, the Aomoto-Drinfel¡¯d-Zagier formula reads \cite{BBBL1997}
\[\sum\limits_{n,m = 1}^\infty  {\zeta \left( {m + 1,{{\left\{ 1 \right\}}_{n - 1}}} \right){x^m}{y^n} = 1 - \exp \left( {\sum\limits_{n = 2}^\infty  {\zeta \left( n \right)\frac{{{x^n} + {y^n} - {{\left( {x + y} \right)}^n}}}{n}} } \right)} ,\tag{4.7}\]
which implies that for any $m,n\in \N$, the multiple zeta value ${\zeta \left( {m + 1,{{\left\{ 1 \right\}}_{n - 1}}} \right)}$ can be
represented as a polynomial of zeta values with rational coefficients, and we have the duality formula
\[\zeta \left( {n + 1,{{\left\{ 1 \right\}}_{m - 1}}} \right) = \zeta \left( {m + 1,{{\left\{ 1 \right\}}_{n - 1}}} \right).\]
In particular, we deduce that
\[\begin{array}{l}
 \zeta \left( {2,{{\left\{ 1 \right\}}_m}} \right) = \zeta \left( {m + 2} \right), \\
 \zeta \left( {3,{{\left\{ 1 \right\}}_m}} \right) = \frac{{m + 2}}{2}\zeta \left( {m + 3} \right) - \frac{1}{2}\sum\limits_{k = 1}^m {\zeta \left( {k + 1} \right)\zeta \left( {m + 2 - k} \right)} .
 \end{array}\]
Therefore, from (4.6) and (4.7), we know that the multiple polylogarithmic values $\zeta \left( {2,{{\left\{ 1 \right\}}_m};\frac{1}{2}} \right)$ can be evaluated in terms of zeta values, polylogarithms and $\ln2$. In particular, one can find explicit formulas for small weights.
\begin{align*}
&\zeta \left( {2;\frac{1}{2}} \right) = \frac{{\zeta \left( 2 \right) - {{\ln }^2}2}}{2},\ \zeta \left( {2,1;\frac{1}{2}} \right) = \frac{1}{8}\zeta \left( 3 \right) - \frac{1}{6}{\ln ^3}2,\\
&\zeta \left( {2,1,1;\frac{1}{2}} \right) = \zeta \left( 4 \right) + \frac{1}{4}\zeta \left( 2 \right){\ln ^2}2 - {\rm{L}}{{\rm{i}}_4}\left( {\frac{1}{2}} \right) - \frac{1}{{12}}{\ln ^4}2 - \frac{7}{8}\ln 2\zeta \left( 3 \right).
\end{align*}
Setting $m=1$ in Theorem 2.2, we have
\[\zeta \left( {\bar 1,\bar 1,{{\left\{ 1 \right\}}_{k - 1}}} \right) =  - {\rm{L}}{{\rm{i}}_{k + 1}}\left( {\frac{1}{2}} \right). \tag{4.8}\]
This formula was also given by \cite{BBBL1997}. Taking $k=1$ in (2.4) and $m=1$ in (2.5), we obtain
\begin{align*}
J\left( {1,m} \right) =& \frac{1}{{m + 1}}{\left( {\ln 2} \right)^{m + 2}} + \sum\limits_{j = 0}^m {{{\left( { - 1} \right)}^{j + 1}}j!{{\left( {\ln 2} \right)}^{m - j}}\left( {\begin{array}{*{20}{c}}
   m  \\
   j  \\
\end{array}} \right)\zeta \left( {j + 2} \right)} \\
&- \sum\limits_{j = 0}^m {\sum\limits_{l = 0}^j {{{\left( { - 1} \right)}^{j + 1}}l!{{\left( {\ln 2} \right)}^{m - l}}\left( {\begin{array}{*{20}{c}}
   m  \\
   j  \\
\end{array}} \right)\left( {\begin{array}{*{20}{c}}
   j  \\
   l  \\
\end{array}} \right){\rm{L}}{{\rm{i}}_{l + 2}}\left( {\frac{1}{2}} \right)} }, \tag{4.9}\\
J\left( {k,1} \right) =& \sum\limits_{i = 0}^k {\sum\limits_{l = 0}^i {{{\left( { - 1} \right)}^i}l!\left( {\begin{array}{*{20}{c}}
   k  \\
   i  \\
\end{array}} \right)\left( {\begin{array}{*{20}{c}}
   i  \\
   l  \\
\end{array}} \right){{\left( {\ln 2} \right)}^{k + 1 - l}}{\rm{L}}{{\rm{i}}_{l + 1}}\left( {\frac{1}{2}} \right)} } \\
& + \sum\limits_{i = 0}^k {\sum\limits_{l = 0}^i {{{\left( { - 1} \right)}^{i + 1}}l!\left( {\begin{array}{*{20}{c}}
   k  \\
   i  \\
\end{array}} \right)\left( {\begin{array}{*{20}{c}}
   i  \\
   l  \\
\end{array}} \right){{\left( {\ln 2} \right)}^{k - l}}\zeta \left( {l + 1,1;\frac{1}{2}} \right)} } .\tag{4.10}
\end{align*}
Hence, putting $k=1$ or $m=2$ in (2.3) and combining (4.9) with (4.10), we get the following corollaries.
\begin{cor} For integer $m\in \N$, we have the recurrence relation
\begin{align*}
\zeta \left( {\bar 1,{{\left\{ 1 \right\}}_{m - 1}},\bar 1} \right) =& \frac{1}{{m!}}{\ln ^{m + 1}}2 - \frac{{\zeta \left( 2 \right){{\ln }^{m - 1}}2}}{{\left( {m - 1} \right)!}} - \sum\limits_{i = 1}^{m - 1} {\frac{{{{\left( {\ln 2} \right)}^i}}}{{i!}}\zeta \left( {\bar 1,{{\left\{ 1 \right\}}_{m - i - 1}},\bar 1} \right)} \\
&- \frac{1}{{\left( {m - 1} \right)!}}\sum\limits_{k = 1}^{m - 1} {\left( {\begin{array}{*{20}{c}}
   {m - 1}  \\
   k  \\
\end{array}} \right){{\left( { - 1} \right)}^{k + 1}}\left\{ \begin{array}{l}
 \sum\limits_{l = 1}^k {l!\left( {\begin{array}{*{20}{c}}
   k  \\
   l  \\
\end{array}} \right){{\left( {\ln 2} \right)}^{m - l - 1}}{\rm{L}}{{\rm{i}}_{l + 2}}\left( {\frac{1}{2}} \right)}  \\
  - k!{\left( {\ln 2} \right)^{m - k - 1}}\zeta \left( {k + 2} \right) \\
 \end{array} \right\}}.\tag{4.11}
\end{align*}
\end{cor}
\begin{cor} For positive integer $k$, then
\begin{align*}
 \zeta \left( {\bar 1,1,\bar 1,{{\left\{ 1 \right\}}_{k - 1}}} \right) =& \ln 2{\rm{L}}{{\rm{i}}_{k + 1}}\left( {\frac{1}{2}} \right) - \frac{{{{\left( { - 1} \right)}^k}}}{{k!}}\left\{ {\sum\limits_{i = 0}^k {\sum\limits_{l = 0}^i {{{\left( { - 1} \right)}^i}l!\left( {\begin{array}{*{20}{c}}
   k  \\
   i  \\
\end{array}} \right)\left( {\begin{array}{*{20}{c}}
   i  \\
   l  \\
\end{array}} \right)} } {{\left( {\ln 2} \right)}^{k + 1 - l}}{\rm{L}}{{\rm{i}}_{l + 1}}\left( {\frac{1}{2}} \right)} \right\} \\
 & - \frac{{{{\left( { - 1} \right)}^k}}}{{k!}}\left\{ {\sum\limits_{i = 0}^k {\sum\limits_{l = 0}^i {{{\left( { - 1} \right)}^{i - 1}}l!\left( {\begin{array}{*{20}{c}}
   k  \\
   i  \\
\end{array}} \right)\left( {\begin{array}{*{20}{c}}
   i  \\
   l  \\
\end{array}} \right)} } {{\left( {\ln 2} \right)}^{k - l}}\zeta \left( {l + 1,1;\frac{1}{2}} \right)} \right\}.\tag{4.12}
\end{align*}
\end{cor}
\section{Proof of Theorem 2.3 and Theorem 2.4}
{\bf Proof of Theorem 2.3}. Similarly as in the proof of Theorem 2.2, we consider the following multiple integral
\begin{align*}
{N_m}\left( k \right):& = \int\limits_0^1 {\frac{1}{{1 + {t_1}}}d{t_1} \cdots \int\limits_0^{{t_{m - 1}}} {\frac{1}{{1 + {t_m}}}d{t_m}\int\limits_0^{{t_m}} {\frac{{{{\ln }^k}\left( {1 + {t_{m + 1}}} \right)}}{{{t_{m + 1}}}}d{t_{m + 1}}} } } \\
 &= \int\limits_{0 < {t_{m + 1}} <  \cdots  < {t_1} < 1} {\frac{{{{\ln }^k}\left( {1 + {t_{m + 1}}} \right)}}{{\left( {1 + {t_1}} \right) \cdots \left( {1 + {t_m}} \right){t_{m + 1}}}}d{t_1} \cdots } d{t_{m + 1}}\\
 &={\left( { - 1} \right)^{k + m}}k!\sum\limits_{{n_1}, \cdots ,{n_{m + 1}} = 1}^\infty  {{{\left( { - 1} \right)}^{{n_1} +  \cdots  + {n_{m + 1}}}}\frac{{{\zeta _{{n_1} - 1}}\left( {{{\left\{ 1 \right\}}_{k - 1}}} \right)}}{{{n_1}}}} \\
 &\quad \quad\quad\quad\quad\quad\times \int\limits_{0 < {t_{m + 1}} <  \cdots  < {t_1} < 1} {t_1^{{n_{m + 1}} - 1} \cdots t_m^{{n_2} - 1}t_{m + 1}^{{n_1} - 1}d{t_1} \cdots } d{t_{m + 1}}\\
 & = {\left( { - 1} \right)^{k + m}}k!\sum\limits_{{n_1}, \cdots ,{n_{m + 1}} = 1}^\infty  {{{\left( { - 1} \right)}^{{n_1} +  \cdots  + {n_{m + 1}}}}\frac{{{\zeta _{{n_1} - 1}}\left( {{{\left\{ 1 \right\}}_{k - 1}}} \right)}}{{n_1^2\left( {{n_1} + {n_2}} \right) \cdots \left( {{n_1} +  \cdots  + {n_{m + 1}}} \right)}}} \\
 & = {\left( { - 1} \right)^{k + m}}k!\zeta \left( {\bar 1,{{\left\{ 1 \right\}}_{m - 1}},2,{{\left\{ 1 \right\}}_{k - 1}}} \right)
 .\tag{5.1}
\end{align*}
By using integration by parts, we can arrive at the conclusion that
\begin{align*}
{N_m}\left( k \right) =& {\left( { - 1} \right)^{m + k - 1}}k!\sum\limits_{i = 1}^{m - 1} {\frac{{{{\left( {\ln 2} \right)}^i}}}{{i!}}\zeta \left( {\bar 1,{{\left\{ 1 \right\}}_{m - i - 1}},2,{{\left\{ 1 \right\}}_{k - 1}}} \right)} \\
& + \frac{{{{\left( { - 1} \right)}^{m - 1}}}}{{\left( {m - 1} \right)!}}\int\limits_0^1 {\frac{{{{\ln }^{m - 1}}\left( {1 + {t_1}} \right)}}{{1 + {t_1}}}d{t_1}} \int\limits_0^{{t_1}} {\frac{{{{\ln }^k}\left( {1 + {t_2}} \right)}}{{{t_2}}}d{t_2}} \\
 = &{\left( { - 1} \right)^{m + k - 1}}k!\sum\limits_{i = 1}^{m - 1} {\frac{{{{\left( {\ln 2} \right)}^i}}}{{i!}}\zeta \left( {\bar 1,{{\left\{ 1 \right\}}_{m - i - 1}},2,{{\left\{ 1 \right\}}_{k - 1}}} \right)} \\
 & + \frac{{{{\left( { - 1} \right)}^{m - 1}}}}{{m!}}{\left( {\ln 2} \right)^m}\int\limits_0^1 {\frac{{{{\ln }^k}\left( {1 + t} \right)}}{t}dt}  + \frac{{{{\left( { - 1} \right)}^m}}}{{m!}}\int\limits_0^1 {\frac{{{{\ln }^{m + k}}\left( {1 + t} \right)}}{t}dt} \\
= &{\left( { - 1} \right)^{m + k - 1}}k!\sum\limits_{i = 1}^{m - 1} {\frac{{{{\left( {\ln 2} \right)}^i}}}{{i!}}\zeta \left( {\bar 1,{{\left\{ 1 \right\}}_{m - i - 1}},2,{{\left\{ 1 \right\}}_{k - 1}}} \right)} \\
 &+ {\left( { - 1} \right)^{m + k - 1}}\frac{{k!}}{{m!}}{\left( {\ln 2} \right)^m}\zeta \left( {\bar 2,{{\left\{ 1 \right\}}_{k - 1}}} \right) + \frac{{{{\left( { - 1} \right)}^k}}}{{m!}}\left( {m + k} \right)!\zeta \left( {\bar 2,{{\left\{ 1 \right\}}_{m + k - 1}}} \right).\tag{5.2}
\end{align*}
Combining (5.1) and (5.2), we obtain (2.6).  This completes the proof of Theorem 2.3.\hfill$\square$\\
{\bf Proof of Theorem 2.4}. By a similar as in the proof of Theorem 2.3, considering the iterated integral
\begin{align*}
{P_m}\left( k \right): =& \int\limits_0^1 {\frac{1}{{1 + {t_1}}}d{t_1} \cdots \int\limits_0^{{t_{m - 2}}} {\frac{1}{{1 + {t_{m - 1}}}}d{t_{m - 1}}\int\limits_0^{{t_{m - 1}}} {\frac{1}{{{t_m}}}d{t_m}} \int\limits_0^{{t_m}} {\frac{{{{\ln }^k}\left( {1 + {t_{m + 1}}} \right)}}{{{t_{m + 1}}}}d{t_{m + 1}}} } } \\
 = &\int\limits_{0 < {t_{m + 1}} <  \cdots  < {t_1} < 1} {\frac{{{{\ln }^k}\left( {1 + {t_{m + 1}}} \right)}}{{\left( {1 + {t_1}} \right) \cdots \left( {1 + {t_{m - 1}}} \right){t_m}{t_{m + 1}}}}d{t_1} \cdots } d{t_{m + 1}}\\
 = &{\left( { - 1} \right)^{k + m + 1}}k!\sum\limits_{{n_1}, \cdots ,{n_m} = 1}^\infty  {{{\left( { - 1} \right)}^{{n_1} +  \cdots  + {n_m}}}\frac{{{\zeta _{{n_1} - 1}}\left( {{{\left\{ 1 \right\}}_{k - 1}}} \right)}}{{n_1^2}}} \\
 &\quad\quad\quad\quad\quad\quad\quad\times \int\limits_{0 < {t_{m + 1}} <  \cdots  < {t_1} < 1} {t_1^{{n_m} - 1} \cdots t_{m - 1}^{{n_2} - 1}t_m^{{n_1} - 1}d{t_1} \cdots } d{t_m}\\
 =& {\left( { - 1} \right)^{k + m + 1}}k!\sum\limits_{{n_1}, \cdots ,{n_m} = 1}^\infty  {{{\left( { - 1} \right)}^{{n_1} +  \cdots  + {n_m}}}\frac{{{\zeta _{{n_1} - 1}}\left( {{{\left\{ 1 \right\}}_{k - 1}}} \right)}}{{n_1^3\left( {{n_1} + {n_2}} \right) \cdots \left( {{n_1} +  \cdots  + {n_m}} \right)}}} \\
 =& {\left( { - 1} \right)^{k + m + 1}}k!\sum\limits_{{n_1} >  \cdots  > {n_m} \ge 1}^\infty  {\frac{{{\zeta _{{n_m} - 1}}\left( {{{\left\{ 1 \right\}}_{k - 1}}} \right)}}{{n_m^3{n_{m - 1}} \cdots {n_1}}}{{\left( { - 1} \right)}^{{n_1}}}} \\
 =& {\left( { - 1} \right)^{k + m + 1}}k!\zeta \left( {\bar 1,{{\left\{ 1 \right\}}_{m - 2}},3,{{\left\{ 1 \right\}}_{k - 1}}} \right).\tag{5.3}
\end{align*}
Similarly, it is easily shown using integration by parts that
\begin{align*}
{P_m}\left( k \right) =& {\left( { - 1} \right)^{k + m}}k!\sum\limits_{i = 1}^{m - 2} {\frac{{{{\left( {\ln 2} \right)}^i}}}{{i!}}\zeta \left( {\bar 1,{{\left\{ 1 \right\}}_{m - i - 2}},3,{{\left\{ 1 \right\}}_{k - 1}}} \right)} \\
&+ \frac{{{{\left( { - 1} \right)}^{m - 2}}}}{{\left( {m - 2} \right)!}}\int\limits_0^1 {\frac{{{{\ln }^{m - 2}}\left( {1 + {t_1}} \right)}}{{1 + {t_1}}}d{t_1}} \int\limits_0^{{t_1}} {\frac{1}{{{t_2}}}d{t_2}\int\limits_0^{{t_2}} {\frac{{{{\ln }^k}\left( {1 + {t_3}} \right)}}{{{t_3}}}d{t_3}} } \\
 = &{\left( { - 1} \right)^{k + m}}k!\sum\limits_{i = 1}^{m - 2} {\frac{{{{\left( {\ln 2} \right)}^i}}}{{i!}}\zeta \left( {\bar 1,{{\left\{ 1 \right\}}_{m - i - 2}},3,{{\left\{ 1 \right\}}_{k - 1}}} \right)} \\
 & + \frac{{{{\left( { - 1} \right)}^{m + k - 2}}}}{{\left( {m - 1} \right)!}}k!{\left( {\ln 2} \right)^{m - 1}}\zeta \left( {\bar 3,{{\left\{ 1 \right\}}_{k - 1}}} \right)\\
 & + \frac{{{{\left( { - 1} \right)}^{m - 1}}}}{{\left( {m - 1} \right)!}}\int\limits_0^1 {\frac{{{{\ln }^{m - 1}}\left( {1 + {t_1}} \right)}}{{{t_1}}}d{t_1}\int\limits_0^{{t_1}} {\frac{{{{\ln }^k}\left( {1 + {t_2}} \right)}}{{{t_2}}}d{t_2}} } .\tag{5.4}
\end{align*}
Hence, by using (5.3) and (5.4), we get
\begin{align*}
\zeta \left( {\bar 1,{{\left\{ 1 \right\}}_{m - 1}},3,{{\left\{ 1 \right\}}_{k - 1}}} \right) =& -\sum\limits_{i = 1}^{m - 1} {\frac{{{{\left( {\ln 2} \right)}^i}}}{{i!}}\zeta \left( {\bar 1,{{\left\{ 1 \right\}}_{m - i - 1}},3,{{\left\{ 1 \right\}}_{k - 1}}} \right)}  - \frac{{{{\left( {\ln 2} \right)}^m}}}{{m!}}\zeta \left( {\bar 3,{{\left\{ 1 \right\}}_{k - 1}}} \right)\\
& + \frac{{{{\left( { - 1} \right)}^k}}}{{m!k!}}\int\limits_0^1 {\int\limits_0^{{t_1}} {\frac{{{{\ln }^m}\left( {1 + {t_1}} \right){{\ln }^k}\left( {1 + {t_2}} \right)}}{{{t_1}{t_2}}}d{t_2}} d{t_1}}.\tag{5.5}
\end{align*}
By using integration by parts again, we find that
\begin{align*}
&\int\limits_0^1 {\int\limits_0^{{t_1}} {\frac{{{{\ln }^m}\left( {1 + {t_1}} \right){{\ln }^k}\left( {1 + {t_2}} \right)}}{{{t_1}{t_2}}}d{t_2}} d{t_1}}  + \int\limits_0^1 {\int\limits_0^{{t_1}} {\frac{{{{\ln }^k}\left( {1 + {t_1}} \right){{\ln }^m}\left( {1 + {t_2}} \right)}}{{{t_1}{t_2}}}d{t_2}} d{t_1}} \\
& = \left( {\int\limits_0^1 {\frac{{{{\ln }^m}\left( {1 + t} \right)}}{t}dt} } \right)\left( {\int\limits_0^1 {\frac{{{{\ln }^k}\left( {1 + t} \right)}}{t}dt} } \right)\\
& = {\left( { - 1} \right)^{m + k}}k!m!\zeta \left( {\bar 2,{{\left\{ 1 \right\}}_{m - 1}}} \right)\zeta \left( {\bar 2,{{\left\{ 1 \right\}}_{k - 1}}} \right).\tag{5.6}
\end{align*}
which together with (5.5) gives (2.7) and finishes the proof of Theorem 2.4. \hfill$\square$\\
Putting $m=k$ in Theorem 2.4, we get the Corollary.
\begin{cor} For positive integer $k$, we have
\begin{align*}
\zeta \left( {\bar 1,{{\left\{ 1 \right\}}_{k - 1}},3,{{\left\{ 1 \right\}}_{k - 1}}} \right) =& -\sum\limits_{i = 1}^{k - 1} {\frac{{{{\left( {\ln 2} \right)}^i}}}{{i!}}\zeta \left( {\bar 1,{{\left\{ 1 \right\}}_{k - i - 1}},3,{{\left\{ 1 \right\}}_{k - 1}}} \right)} \\
& - \frac{{{{\left( {\ln 2} \right)}^k}}}{{k!}}\zeta \left( {\bar 3,{{\left\{ 1 \right\}}_{k - 1}}} \right) + {\left( { - 1} \right)^k}\frac{1}{2}{\zeta ^2}\left( {\bar 2,{{\left\{ 1 \right\}}_{k - 1}}} \right).\tag{5.7}
\end{align*}
\end{cor}
\section{Some results and cases}
From Corollary 3.3 and 4.1, we get the following Theorem.
\begin{thm} For integers $m,k\in \N$, then the multiple zeta values \[\zeta (\bar 1,{\left\{ 1 \right\}_{m - 1}},\bar 1)\]and \[\zeta (\bar 1,{\left\{ 1 \right\}_{m - 1}},2,{\left\{ 1 \right\}_{k - 1}})\]
are expressible in terms of zeta values, polylogarithms and $\ln2$.
\end{thm}
From Theorem 2.1-2.4, we can give the following results.
\begin{ex}
Some closed form of integrals $I(k,m)$ and $J(k,m)$
\begin{align*}
 &J(0,1) = \frac{1}{2}{\ln ^2}2, \\
 &I(1,1) =  - \frac{3}{4}\zeta \left( 3 \right), \\
 &J(1,1) = \frac{1}{3}{\ln ^3}2 - \frac{1}{2}\zeta \left( 2 \right)\ln 2 + \frac{1}{8}\zeta \left( 3 \right) \\
 &J(2,1)= \frac{1}{4}{\ln ^4}2 + 2\zeta \left( 3 \right)\ln 2 - \zeta \left( 2 \right){\ln ^2}2 - \frac{1}{4}\zeta \left( 4 \right) \\
 &I(0,3) = 6\zeta \left( 4 \right) + \frac{3}{2}\zeta \left( 2 \right){\ln ^2}2 - \frac{1}{4}{\ln ^4}2 - \frac{{21}}{4}\zeta \left( 3 \right)\ln 2 - 6{\rm{L}}{{\rm{i}}_4}\left( {\frac{1}{2}} \right), \\
 &I(1,2)= \frac{{15}}{4}\zeta \left( 4 \right) + \zeta \left( 2 \right){\ln ^2}2 - \frac{1}{6}{\ln ^4}2 - \frac{7}{2}\zeta \left( 3 \right)\ln 2 - 4{\rm{L}}{{\rm{i}}_4}\left( {\frac{1}{2}} \right), \\
 &J(1,2) = \frac{1}{3}{\ln ^4}2{\rm{ + }}2\zeta \left( 3 \right)\ln 2{\rm{ + }}2{\rm{L}}{{\rm{i}}_4}\left( {\frac{1}{2}} \right) - \zeta \left( 2 \right){\ln ^2}2 - 2\zeta \left( 4 \right), \\
 &I(0,4) =  - 24{\rm{L}}{{\rm{i}}_5}\left( {\frac{1}{2}} \right) - 24\ln 2{\rm{L}}{{\rm{i}}_4}\left( {\frac{1}{2}} \right) - \frac{4}{5}{\ln ^5}2 - \frac{{21}}{2}\zeta \left( 3 \right){\ln ^2}2 + 24\zeta \left( 5 \right) + 4\zeta \left( 2 \right){\ln ^3}2.
\end{align*}
\end{ex}
\begin{ex}
Some closed form of multiple zeta values
\begin{align*}
&\zeta \left( {\bar 1,1,\bar 1} \right) = \frac{1}{8}\zeta \left( 3 \right) - \frac{1}{6}{\ln ^3}2,\\
&\zeta \left( {\bar 1,2} \right) = \frac{1}{2}\zeta \left( 2 \right)\ln 2 - \frac{1}{4}\zeta \left( 3 \right),\\
&\zeta \left( {\bar 1,3} \right) = \frac{3}{4}\zeta \left( 3 \right)\ln 2 - \frac{5}{{16}}\zeta \left( 4 \right),\\
&\zeta \left( {\bar 2,1,1} \right) = {\rm{L}}{{\rm{i}}_4}\left( {\frac{1}{2}} \right) + \frac{1}{{24}}{\ln ^4}2 + \frac{7}{8}\zeta \left( 3 \right)\ln 2 - \frac{1}{4}\zeta \left( 2 \right){\ln ^2}2 - \zeta \left( 4 \right),\\
&\zeta \left( {\bar 1,1,2} \right) = 3{\rm{L}}{{\rm{i}}_4}\left( {\frac{1}{2}} \right) + \frac{1}{8}{\ln ^4}2 + \frac{{23}}{8}\zeta \left( 3 \right)\ln 2 - \zeta \left( 2 \right){\ln ^2}2 - 3\zeta \left( 4 \right),\\
&\zeta \left( {\bar 1,1,1,\bar 1} \right) = {\rm{L}}{{\rm{i}}_4}\left( {\frac{1}{2}} \right) + \frac{1}{{12}}{\ln ^4}2 + \frac{7}{8}\zeta \left( 3 \right)\ln 2 - \frac{1}{2}\zeta \left( 2 \right){\ln ^2}2 - \zeta \left( 4 \right).
\end{align*}
\end{ex}
Next, we now close this paper with two Theorem.
\begin{thm} For integers $m\geq 1$ and $k\geq 0$, then the following identity holds:
\begin{align*}
&m!\sum\limits_{l = 0}^k {l!\left( {\begin{array}{*{20}{c}}
   k  \\
   l  \\
\end{array}} \right){{\left( {\ln 2} \right)}^{k - l}}\zeta \left( {l + 2,{{\left\{ 1 \right\}}_{m - 1}};\frac{1}{2}} \right)} \\
& + k!\sum\limits_{l = 0}^m {l!\left( {\begin{array}{*{20}{c}}
   m  \\
   l  \\
\end{array}} \right){{\left( {\ln 2} \right)}^{m - l}}\zeta \left( {l + 1,{{\left\{ 1 \right\}}_{k }};\frac{1}{2}} \right)} \\
&= m!k!\zeta \left( {m + 1,{{\left\{ 1 \right\}}_k}} \right).\tag{6.1}
\end{align*}
\end{thm}
\pf By using (3.1), (3.6) and (3.7), we deduce that
\begin{align*}
\int\limits_0^{1/2} {\frac{{{{\ln }^k}x{{\ln }^m}\left( {1 - x} \right)}}{x}dx}  =& {\left( { - 1} \right)^m}m!\sum\limits_{n = 1}^\infty  {\frac{{{\zeta _{n - 1}}\left( {{{\left\{ 1 \right\}}_{m - 1}}} \right)}}{n}\int\limits_0^1 {{x^{n - 1}}{{\ln }^k}xdx} } \\
= &{\left( { - 1} \right)^{m + k}}m!\sum\limits_{l = 0}^k {l!\left( {\begin{array}{*{20}{c}}
   k  \\
   l  \\
\end{array}} \right){{\left( {\ln 2} \right)}^{k - l}}\zeta \left( {l + 2,{{\left\{ 1 \right\}}_{k - 1}};\frac{1}{2}} \right)} \\
&\mathop  = \limits^{1 - x = t} \int\limits_{1/2}^1 {\frac{{{{\ln }^k}\left( {1 - t} \right){{\ln }^m}t}}{{1 - t}}dt}\\
 =& {\left( { - 1} \right)^k}k!\sum\limits_{n = 1}^\infty  {{\zeta _{n - 1}}\left( {{{\left\{ 1 \right\}}_k}} \right)\int\limits_{1/2}^1 {{t^{n - 1}}{{\ln }^m}tdt} } \\
 =& {\left( { - 1} \right)^{m + k}}m!k!\zeta \left( {m + 1,{{\left\{ 1 \right\}}_k}} \right)\\& + {\left( { - 1} \right)^{m + k - 1}}k!\sum\limits_{l = 0}^m {l!\left( {\begin{array}{*{20}{c}}
   m  \\
   l  \\
\end{array}} \right){{\left( {\ln 2} \right)}^{m - l}}\zeta \left( {l + 1,{{\left\{ 1 \right\}}_{k }};\frac{1}{2}} \right)}.\tag{6.2}
\end{align*}
By simple calculation, we obtain formula (6.1).\hfill$\square$\\
Taking $m=2,k=1$ in (6.2), we have
\[\zeta \left( {3,1;\frac{1}{2}} \right) = \frac{1}{8}\zeta \left( 4 \right) - \frac{1}{8}\zeta \left( 3 \right)\ln 2 + \frac{1}{{24}}{\ln ^4}2.\]
Proceeding in a similar fashion to evaluation of the Theorem 2.1-2.4, it is possible to evaluate other alternating multiple zeta values. For example, by using (3.1) and applying the same arguments as in the proof of Theorem 2.1 and 2.2, we may easily deduce the following integral representation of series
\[\zeta \left( {\overline {m + 1},\bar 1,{{\left\{ 1 \right\}}_{k - 1}}} \right) = \frac{{{{\left( { - 1} \right)}^{m + k - 1}}}}{{k!m!}}\int\limits_0^1 {\frac{{{{\left( {\ln t} \right)}^m}{{\ln }^k}\left( {1 - t} \right)}}{{1 + t}}dt},\tag{6.3} \]
and identities
\begin{align*}
\zeta \left( {\bar 1,{{\left\{ 1 \right\}}_{m - 1}},\bar 1,\bar 1,{{\left\{ 1 \right\}}_{k - 1}}} \right) =& \frac{{{{\left( { - 1} \right)}^{k - 1}}}}{{m!k!}}\left\{ {k{{\left( {\ln 2} \right)}^m}J\left( {1,k - 1} \right) - \left( {m + k} \right)J\left( {1,m + k - 1} \right)} \right\}\\
& - \sum\limits_{i = 1}^{m - 1} {\frac{{{{\left( {\ln 2} \right)}^i}}}{{i!}}} \zeta \left( {\bar 1,{{\left\{ 1 \right\}}_{m - i - 1}},\bar 1,\bar 1,{{\left\{ 1 \right\}}_{k - 1}}} \right),\tag{6.4}\\
\zeta \left( {\bar 1,{{\left\{ 1 \right\}}_{m - 1}},\bar 1,\bar 1,\bar 1,{{\left\{ 1 \right\}}_{k - 1}}} \right)& = {\left( { - 1} \right)^{m + 1}}\zeta \left( {k + 1,2,{{\left\{ 1 \right\}}_{m - 1}};\frac{1}{2}} \right).\tag{6.5}
\end{align*}
Hence, from (4.9) and (6.4), we obtain the conclusion.
\begin{thm}
 If $m,k\in \N$, then the alternating multiple zeta values $\zeta \left( {\bar 1,{{\left\{ 1 \right\}}_{m - 1}},\bar 1,\bar 1,{{\left\{ 1 \right\}}_{k - 1}}} \right)$
can be expressed as a rational linear combination of zeta values, polylogarithms and $\ln2$.
\end{thm}
From (6.4) and (6.5), we can get the following results
\begin{align*}
&\zeta \left( {\bar 1,\bar 1,\bar 1,\bar 1} \right) = \frac{1}{{24}}{\ln ^4}2 + \frac{1}{4}\zeta \left( 3 \right)\ln 2 - \frac{1}{4}\zeta \left( 2 \right){\ln ^2}2 + \frac{1}{{16}}\zeta \left( 4 \right),\\
&\zeta \left( {\bar 1,\bar 1,\bar 1,1} \right) = 3{\rm{L}}{{\rm{i}}_4}\left( {\frac{1}{2}} \right) + \frac{1}{6}{\ln ^4}2 + \frac{{23}}{8}\zeta \left( 3 \right)\ln 2 - \zeta \left( 2 \right){\ln ^2}2 - 3\zeta \left( 4 \right),\\
&\zeta \left( {\bar 1,1,\bar 1,\bar 1} \right) =  - 3{\rm{L}}{{\rm{i}}_4}\left( {\frac{1}{2}} \right) - \frac{1}{{12}}{\ln ^4}2 - \frac{{11}}{4}\zeta \left( 3 \right)\ln 2 - \frac{3}{4}\zeta \left( 2 \right){\ln ^2}2 + 3\zeta \left( 4 \right).
\end{align*}
{\bf Acknowledgments.} The authors would like to thank the anonymous
referee for his/her helpful comments, which improve the presentation
of the paper.
 \end{document}